\documentclass[12pt]{amsart}
\usepackage{latexsym}
\usepackage{a4}
\usepackage{amssymb}
\usepackage{amsbsy}
\usepackage{bm}
\newtheorem{thm}{Theorem}[section] 
\newtheorem{pro}[thm]{Proposition} 
\newtheorem{lem}[thm]{Lemma} 
 
\newtheorem{cor}[thm]{Corollary} 
\theoremstyle{definition}
 
\newtheorem{de}[thm]{Definition}

\numberwithin{equation}{section} 

\newcommand{\ab}[1]{{\mathbf{#1}}}

\newcommand{\N}{\Bbb{ N}}

\newcommand{\setsuchthat}{\,\, \pmb{|} \,\,}

\newcommand{\R}{\mathbb{R}}

\newcommand{\tup}[2]{#1_{1},\ldots,#1_{#2}}
\DeclareMathAlphabet{\mathbfsl}{OT1}{cmr}{bx}{it}
\newcommand{\tupBold}[1]{\mathbfsl{#1}}
\newcommand{\vb}[1]{\tupBold{#1}} 

\newcommand{\rad}{\mathrm{rad}}
\newcommand{\radchi}{\rad_{\chi}}
\newcommand{\clo}[1]{\overline{#1}}
\newcommand{\FGen}[1]{({#1})}
\newcommand{\RGen}[1]{[{#1}]}
\newcommand{\ran}{\mathrm{range}}
\newcommand{\Irr}{\mathrm{Irr}}
\renewcommand{\emptyset}{\varnothing}
\title{On function compositions that are polynomials}
\author{Erhard Aichinger}
\address{Erhard Aichinger,
Institut f\"ur Algebra,
Johannes Kepler Universit\"at Linz,
4040 Linz,
Austria}
\email{erhard@algebra.uni-linz.ac.at}
\urladdr{http://www.algebra.uni-linz.ac.at}
\thanks{Supported by the Austrian Science Fund (FWF):P24077}
\begin{document}
\bibliographystyle{amsplain}
\begin{abstract}
   For a polynomial map $\tupBold{f} : k^n \to k^m$ ($k$ a field),
   we investigate those polynomials $g \in k[t_1,\ldots, t_n]$
   that can be written as a composition $g = h \circ \tupBold{f}$,
   where $h: k^m \to k$ is an arbitrary function.
   In the case that $k$ is algebraically closed of characteristic $0$ and
   $\tupBold{f}$ is surjective, we will show that 
   $g = h \circ \tupBold{f}$ implies that $h$ is 
   a polynomial.
 \end{abstract} 
\maketitle
 \section{Introduction}
  In the present note we investigate the situation that the
  value of a polynomial depends only on the value of certain
  given polynomials. To be precise,
  let $k$ be a field, $m, n \in \N$, and let
 $g, f_1, \ldots, f_m \in k[\tup{t}{n}]$.
 We say that $g$ \emph{is determined by $\tupBold{f} = (\tup{f}{m})$} if
 for all $\tupBold{a}, \tupBold{b} \in k^n$ with
 $f_1 (\tupBold{a}) = f_1 (\tupBold{b})$, \ldots,
 $f_m (\tupBold{a}) = f_m (\tupBold{b})$, we have
 $g (\tupBold{a}) = g(\tupBold{b})$. In other words,
 $g$ is determined by $\tupBold{f}$ if and only if there is
 a function $h : k^m \to k$ such that
 \[
     g (\tupBold{a}) = h (f_1 (\tupBold{a}), \ldots, f_m (\tupBold{a}))  \text{ for all } \tupBold{a} \in k^n.
 \]
 For given $f_1,\ldots, f_m \in k[\tup{t}{n}]$, the set
 of all elements of $k[\tup{t}{n}]$ that are determined by $(f_1,\ldots,f_m)$ is
 a $k$-subalgebra of $k[\tup{t}{n}]$; we will denote this $k$-subalgebra
 by $k \langle \tup{f}{m} \rangle$ or $k \langle \tupBold{f} \rangle$.
 As an example, we see that $t_1 \in \R \langle {t_1}^3 \rangle$; more generally,
 if $(f_1,\ldots, f_m) \in k[t_1,\ldots, t_n]^m$ induces an injective map
 from $k^n$ to $k^m$, we have $k \langle \tupBold{f} \rangle = k [t_1,\ldots, t_n]$. 
 In the present note, we will describe $k \langle \tupBold{f} \rangle$
   in the case that $k$
  is algebraically closed and $\tupBold{f}$ induces a map from $k^n$ to $k^m$
  that is surjective, or, in a sense specified later, at least close to being
  surjective.

 The first set that $k \langle \tupBold{f} \rangle$ is compared with is    
  the $k$-subalgebra of $k[\tup{t}{n}]$ generated  by $\{\tup{f}{m}\}$,
  which we will denote by $k \RGen{\tup{f}{m}}$ or $k \RGen{\tupBold{f}}$;
 in this algebra, we find exactly those polynomials that can be
  written as $p(f_1,\ldots, f_m)$ with $p \in k[x_1,\ldots, x_m]$.
  Clearly, $k \RGen{\tupBold{f}} \subseteq k \langle \tupBold{f} \rangle$.
   The other inclusion need not hold in general: on any field $k$,
 let $f_1 = t_1$, $f_2 = t_1 t_2$. Then $\frac{f_2^2}{f_1} = t_1 t_2^2$ is
 $(f_1,f_2)$-determined, but $t_1 t_2^2 \not\in k \RGen{f_1,f_2}$.  

  The second set with which we will compare $k \langle \tupBold{f} \rangle$
  is the set of all polynomials that can be written as rational functions
  in $f_1,\ldots, f_m$.
  We denote the quotient field
     of $k[t_1,\ldots, t_n]$ by $k (t_1,\ldots, t_n)$. For  $r_1,\ldots, r_m \in 
     k(t_1,\ldots,t_n)$, the subfield of $k(t_1,\ldots, t_n)$ that is 
     generated by $k \cup \{r_1,\ldots, r_m\}$ is denoted
     $k \FGen{r_1,\ldots,r_m}$.
  We first observe that there are polynomials that can be written as rational
  functions in $\tupBold{f}$, but fail to be $\tupBold{f}$-determined.
  As an example, we see that $t_2 \in k \FGen{t_1, t_1 t_2}$, but since
     $(0, 0 \cdot 0) = (0, 0 \cdot 1)$ and $0 \neq 1$, the polynomial
     $t_2$ is not $(t_1,  t_1 t_2)$-determined. 
  As for the converse inclusion, we take a field $k$ of positive characteristic $\chi$.
  Then $t_1$ is $({t_1}^{\chi})$-determined, but $t_1 \not\in k \FGen{{t_1}^{\chi}}$.

 On the positive side, it is known that
    $k\RGen{\tup{f}{m}} = k \langle \tup{f}{m} \rangle$
 holds in the following cases (cf. \cite[Theorem~3.1]{AS:APOA}):
 \begin{itemize}
    \item $k$ is algebraically closed, $m=n=1$, and the derivative
          $f'$ of $f$ is not the zero polynomial, and, more generally,
    \item $k$ is algebraically closed, $m = n$, and there are
          univariate polynomials $g_1,\ldots, g_m \in k[t]$ with
          $g_1' \neq 0, \ldots, g_m' \neq 0$, $f_1 = g_1 (t_1), \ldots,
           f_m = g_m (t_m)$.
  \end{itemize}
  Let us now briefly outline the results obtained in the present note:
  Let  $k$ be an algebraically closed field of characteristic~$0$, and
  let $f_1,\ldots, f_m \in k[t_1,\ldots, t_n]$ be algebraically independent over $k$.
   Then we have
  $k \langle \tupBold{f} \rangle \subseteq k \FGen{ \tupBold{f} }$ 
  (Theorem~\ref{thm:t1}).
  The equality $k \RGen{ \tupBold{f} } = k \langle \tupBold{f} \rangle$
  holds if and only if $\tupBold{f}$ induces a map from $k^n$ to $k^m$
  that is \emph{almost surjective} (see Definition~\ref{de:as}).
  This equality is stated in Theorem~\ref{thm:isinring0}.
  Similar results are given for the case of positive characteristic.

   The last equality has a consequence on the functional decomposition 
   of polynomials. If
$\tupBold{f}$ induces a surjective mapping from $k^n$ to $k^m$, ($k$ algebraically
  closed of characteristic $0$), and if 
$h :k^m \to k$ is an arbitrary function such that 
$h \circ \tupBold{f}$ is a polynomial function, then $h$ is a polynomial function.
In an algebraically closed field of positive characteristic $\chi$, we will
conclude that $h$ is a composition of taking $\chi$-th roots and a 
polynomial function (Corollary~\ref{cor:comp}).

\section{Preliminaries about polynomials}
    For the notions from algebraic geometry used in this note,
    we refer to \cite{CLO:IVA}; deviating from their definitions,
    we call the set of solutions of a system of polynomial equations
    an \emph{algebraic set} (instead of \emph{affine variety}). For an
    algebraically closed field $k$ and 
    $A \subseteq k^m$, we let $I_m (A)$ (or simply $I(A)$) be the set
    of polynomials vanishing on every point in $A$, and for $P \subseteq k[t_1,\ldots, t_m]$,
    we let $V_m (P)$ (or simply $V (P)$) be the set of common zeroes of $P$ in $k^m$.
    The Zariski-closure $V(I(A))$ of a set $A \subseteq k^m$ will be abbreviated
    by $\overline{A}$. The \emph{dimension} of an algebraic set $A$ is
    the maximal $d \in \{0, \ldots, m\}$ such that there are $i_1 < i_2 < \cdots < i_d \in
    \{1,\ldots, m\}$ with 
    $I(A) \cap k [x_{i_1}, \ldots, x_{i_d}] = \{0\}$. We abbreviate the dimension of
    $A$ by $\dim (A)$ and set $\dim (\emptyset) := -1$.
     For $f_1, \ldots, f_m, g \in k[t_1, \ldots, t_n]$, and 
     $D := \{ (f_1 (\vb{a}), \ldots, f_m (\vb{a}), g(\vb{a})) \setsuchthat \vb{a} \in k^n\}$,
     its Zariski-closure $\clo{D}$ is an irreducible algebraic set, and its dimension
     is the maximal number of algebraically independent elements in
     $\{f_1, \ldots, f_m, g\}$. The Closure Theorem \cite[p.\ 258]{CLO:IVA} tells
     that there exists an algebraic set $W \subseteq k^{m+1}$ with
     $\dim (W) < \dim (\clo{D})$ such that $\clo{D} = D \cup W$.
     If $\dim (\clo{D}) = m$, then there exists an irreducible polynomial
     $p \in k[x_1, \ldots, x_{m+1}]$ such that
     $\clo{D} = V(p)$. We will denote this $p$ by $\Irr (\clo{D})$; $\Irr (\clo{D})$
     is then defined up to a multiplication with a nonzero element from $k$.
     Above this, we recall that 
    a set is \emph{constructible} if and only if it can be generated
    from algebraic sets by a finite application of the set-theoretic operations of
    forming the union of two sets, the intersection of two sets, and the complement of a set,
      and that the range of a polynomial map from
    $k^n$ to $k^m$ and its complement are constructible. This
    is of course a consequence of the Theorem of Chevalley-Tarski
    \cite[Exercise~II.3.19]{Ha:AG}, but since we are only concerned
    with the image of $k^n$, it also follows from \cite[p.262, Corollary~2]{CLO:IVA}. 
    \begin{de} \label{de:as}
       Let $k$ be an algebraically closed field, $m,n \in \N$, and
       let $\tupBold{f} = (f_1,\ldots, f_m) \in (k[t_1,\ldots, t_n])^m$.   
       By $\ran (\tupBold{f})$, we denote the image
       of the mapping $\hat{\tupBold{f}} : k^n \to k^m$ that
       is induced by $\tupBold{f}$.
       We say that $\tupBold{f}$ is \emph{almost surjective on $k$}   
       if the dimension of the Zariski-closure of
       $k^m \setminus \ran (\tupBold{f})$ is at most
       $m-2$.
    \end{de}
  \begin{pro} \label{pro:asai}
    Let $k$ be an algebraically closed field, and 
    let $(f_1,\ldots, f_m) \in k[t_1,\ldots,t_n]^m$ be  almost surjective on $k$.
   Then the sequence
  $(f_1, \ldots, f_m)$ is algebraically independent over $k$.
   \end{pro}
   \emph{Proof:}
    Seeking
    a contradiction, we suppose that there is $u \in k[x_1,\ldots, x_m]$ with
    $u \neq 0$ and $u (f_1,\ldots, f_m) = 0$. Then
    $\ran (\tupBold{f}) \subseteq V(u)$, hence
    $\dim (\clo{\ran (\tupBold{f})}) \le m-1$.
    Since $\tupBold{f}$ is almost surjective, $k^m$ is then the
    union of two algebraic sets of dimension $\le m-1$, a contradiction.  \qed

  We will use the following easy consequence of the description of
  constructible sets:
    \begin{pro} \label{pro:BWX}
      Let $k$ be an algebraically closed field, and let $B$ be a constructible subset of $k^m$
      with $\dim (\clo{B}) \ge m-1$. Then  there exist algebraic sets
     $W, X$ such that $W$ is irreducible, $\dim (W) = m-1$, $\dim (X) \le m-2$, 
      $W \setminus X \subseteq B$.
    \end{pro}
      Since $B$ is constructible, there
      are irreducible algebraic sets $V_1, \ldots, V_p$ and algebraic sets
      $W_1, \ldots, W_p$ with $W_i \subsetneq V_i$ and $B = \bigcup_{i=1}^p (V_i \setminus W_i)$
       (cf. \cite[p.\ 262]{CLO:IVA}).
      We assume that the $V_i$'s are ordered with nonincreasing dimension.
      If $\dim (V_1) = m$, then $k^m \setminus W_1 \subseteq B$.
      Let $U$ be an irreducible algebraic set of dimension $m-1$ with
      $U \not\subseteq W_1$.
      Then $U \cap (k^m \setminus W_1) =
            U \setminus (W_1 \cap U)$. Since $W_1 \cap U \neq U$, setting
      $W := U, X := W_1 \cap U$ yields the required sets.
   
      If $\dim (V_1) = m-1$, then $W := V_1$ and $X := W_1$ are the required sets.
      
      The case $\dim (V_1) \le m-2$ cannot occur because
      then $\clo{B} \subseteq V_1 \cup \ldots \cup V_p$ has dimension at most $m-2$. \qed
  
     Let $k$ be a field, and let $p,q, f \in k[t]$ such that
     $\deg (f) > 0$. It is known that $p(f)$ divides $q(f)$ if and
    only if $p$ divides $q$ \cite[Lemma~2.1 and 2.2]{En:PS}. The following Lemma
     yields a multivariate version of this result. 
    \begin{lem} \label{lem:polies}
      Let $k$ be an algebraically closed field, $m, n \in \N$, and
      let $\tupBold{f} = (f_1,\ldots, f_m)$ $\in (k[t_1,\ldots, t_n])^m$.   
      Then the following are equivalent:
      \begin{enumerate}
         \item \label{it:l0} $\tupBold{f}$ is almost surjective on $k$.
         \item \label{it:l1} $k \FGen{f_1,\ldots, f_m} \cap k[t_1,\ldots, t_n]
                = k \RGen{f_1,\ldots, f_m}$ and $(f_1,\ldots, f_m)$ is algebraically
                independent over $k$.
         \item \label{it:l2}
                For all $p,q \in k[x_1,\dots, x_m]$ with
                $p(f_1,\ldots, f_m) \mid q (f_1, \ldots, f_m)$, we have
                $p \mid q$.
      \end{enumerate}
    \end{lem}
     \emph{Proof:} 
        \eqref{it:l0}$\Rightarrow$\eqref{it:l1}:
        (This proof uses some ideas from the proof of Theorem~4.2.1 in
        \cite[p.\ 82]{vdE:PAAT}.)
         Let $g \in k \FGen{f_1,\ldots,f_m} \cap k[t_1,\ldots, t_n]$.
        Then there are
        $r,s \in k[x_1,\ldots,x_m]$ with
   $\gcd (r,s) = 1$ and $g = r(f_1,\ldots,f_m) / s (f_1,\ldots, f_m)$, and thus
   \begin{equation} \label{eq:g}
      g (t_1,\ldots, t_n) \cdot s(f_1,\ldots,f_m) =
                                r (f_1,\ldots, f_m).
   \end{equation}
    Suppose $s \not\in k$. Then
   $V(s)$ has dimension $m-1$.
   We have
   $V(s) = (V(s) \cap \ran (\tupBold{f}))
           \cup (V(s) \cap (k^m \setminus \ran (\tupBold{f})))
         \subseteq \clo{V(s) \cap \ran (\tupBold{f})}
                   \cup
                   \clo{V(s) \cap (k^m \setminus \ran (\tupBold{f}))}$.
  Since $\tupBold{f}$ is almost surjective, 
   $\clo{V(s) \cap \ran (\tupBold{f})}$ 
   is then of dimension $m-1$.
   Hence it contains an irreducible component of dimension $m-1$,
   and thus there is an irreducible $p \in k[x_1,\ldots, x_m]$
   such that $V(p) \subseteq \clo{V(s) \cap \ran (\tupBold{f})}$.
   Since then $V(p) \subseteq V(s)$, the Nullstellensatz yields $n_1 \in \N$ 
   with $p \mid s^{n_1}$, and thus by the irreducibility of $p$, $p \mid s$.
   Now we show that for all $\vb{a} \in V(s) \cap \ran (\tupBold{f})$, we have
   $r (\vb{a}) = 0$. To this end, let $\vb{b} \in k^n$ with
   $\tupBold{f} (\vb{b}) = \vb{a}$.
   Setting $\vb{t} := \vb{b}$ in~\eqref{eq:g}, we obtain
   $r (\vb{a}) = 0$.
   Thus $V(s) \cap \ran (\tupBold{f}) \subseteq V(r)$, and therefore
   $\clo{V(s) \cap \ran (\tupBold{f})} \subseteq V(r)$, which implies
   $V(p) \subseteq V(r)$. By the Nullstellensatz,
   we have an $n_2 \in \N$ with $p \mid r^{n_2}$ and thus
   by the irreducibility of $p$, $p \mid r$.
   Now $p \mid r$ and $p \mid s$, contradicting $\gcd (r,s) = 1$.   
   Hence $s \in k$,  and thus $g \in k \RGen{f_1,\ldots, f_m}$. 
   The algebraic independence of $(f_1,\ldots, f_m)$ follows from
   Proposition~\ref{pro:asai}.  

   \eqref{it:l1}$\Rightarrow$\eqref{it:l2}: 
   Let $p,  q \in k[x_1,\ldots, x_m]$ such that
    $p  (f_1,\ldots, f_m) \mid q (f_1,\ldots, f_m)$.
    If $p (f_1,\ldots, f_m) = 0$, then 
    $q (f_1,\ldots, f_m) = 0$, and thus by the algebraic independence
    of $(f_1,\ldots, f_m)$, we have $q=0$ and thus $p \mid q$.
    Now assume $p(f_1,\ldots, f_m) \neq 0$. We have
    $a (t_1, \ldots, t_n) \in k[t_1,\ldots, t_n]$ such that
      \begin{equation} \label{eq:qa} 
       q (f_1,\ldots, f_m) = a (t_1,\ldots,t_n) \cdot p (f_1,\ldots, f_m) 
      \end{equation}
     and thus $a(t_1,\ldots, t_n) \in k \FGen{f_1,\ldots, f_m} \cap 
                                      k [t_1,\ldots, t_n]$.
     Thus there exists $b \in k[x_1,\ldots, x_m]$ such that
     $a (t_1,\ldots, t_n) = b (f_1,\ldots, f_m)$.
     Now \eqref{eq:qa} yields
       \[ q (f_1,\ldots, f_m) = b (f_1,\ldots, f_m) \cdot  p (f_1,\ldots, f_m). \]
     Using the algebraic independence of $(f_1,\ldots, f_m)$, we obtain
     $q(x_1,\ldots, x_m) = b(x_1, \ldots, x_m) \cdot p (x_1, \ldots, x_m)$, 
     and thus $p \mid q$.

    \eqref{it:l2}$\Rightarrow$\eqref{it:l0}
        Seeking a contradiction, we suppose that $\tupBold{f}$ is not almost surjective. 
      Let $B := k^m \setminus \ran (\tupBold{f})$. Then
      $\dim (\clo{B}) \ge m-1$. Since $B$ is constructible, Proposition~\ref{pro:BWX}
      yields
      $W, X$ with $W$ irreducible,  $\dim (W) = m-1$, $\dim (X) \le m-2$, 
      $W \setminus X \subseteq B$. Since $W$ is irreducible of dimension $m-1$, there is
      $p \in k[x_1,\ldots, x_m]$ such that $W = V(p)$.
      Since $\dim (W) > \dim (X)$, we have $W \not\subseteq X$,
      thus $I(X) \not\subseteq I(W)$, and therefore there is
      $q \in I(X)$ with $q \not\in I(W)$.
      We have $W \subseteq B \cup X$, and thus
      $W \cap \ran (\tupBold{f}) \subseteq X$.
      This implies that for all $\tupBold{a} \in k^n$ with
      $p(\tupBold{f} (\tupBold{a})) = 0$, we have
      $q(\tupBold{f} (\tupBold{a})) = 0$: in fact, 
      if $p (\tupBold{f} (\tupBold{a})) = 0$, then
      $\tupBold{f} (\ab{a}) \in V(p) \cap
                       \ran (\tupBold{f}) = 
               W \cap \ran (\tupBold{f}) \subseteq X$.
     Hence $q(\tupBold{f} (\tupBold{a})) = 0$.
     By the Nullstellensatz, we obtain a $\nu \in \N$ such that
     $p (f_1,\ldots, f_m) \mid q(f_1,\ldots, f_m)^\nu$.
     Therefore, using~\eqref{it:l2}, we have $p \mid q^\nu$.     
     This implies $V(p) \subseteq V(q)$. Thus we have
     $W \subseteq V(q)$ and therefore $q \in I(W)$, contradicting
     the choice of $q$.
     Hence $\tupBold{f}$ is almost surjective, proving~\eqref{it:l0}. \qed

 \section{$\tupBold{f}$-determined polynomials} \label{sec:fdet}

 We will first show that often all $\tupBold{f}$-determined polynomials
 are rational functions of $\tupBold{f}$. Special care, however, is needed in the 
 case of positive characteristic. In an algebraically closed field of characteristic
 $\chi > 0$, the unary polynomial $t_1$ is $(t_1^\chi)$-determined, but $t_1$
 is neither a polynomial nor a rational function of $t_1^{\chi}$. 
 \begin{de}
    Let $k$ be a field of characteristic $\chi > 0$, let $n \in \N$, and
    let $P$ be a subset of $k[t_1,\ldots, t_n]$.
    We define the set $\radchi (P)$ by
    \[
         \radchi (P) := \{f \in k[t_1,\ldots, t_n] \setsuchthat 
                          \text{ there is $\nu \in \N_0$ such that } f^{\chi^\nu} \in P\}.
    \]
 \end{de}
 \begin{lem} \label{lem:m}
      Let $k$ be an algebraically closed field, let $m, n \in \N$, let 
      $f_1,\ldots, f_m$ be algebraically
      independent polynomials in $k[t_1,\ldots,t_n]$,  let
      $g \in k \langle f_1,\ldots, f_m \rangle$, and let
      $D := \{ (f_1 (\vb{a}), \ldots, f_m (\vb{a}), g (\vb{a})) \setsuchthat
               \vb{a} \in k^n \}$.
      Then $\dim (\clo{D}) = m$.
\end{lem}
 \emph{Proof:}
    By the closure theorem \cite[p.\ 258]{CLO:IVA}, there is an algebraic set $W$ such that
    $\clo{D} = D \cup W$ and $\dim (W) < \dim (\clo{D})$.
    Let $\pi : k^{m+1} \to k^m, (y_1,\ldots, y_{m+1}) \mapsto
                               (y_1,\ldots, y_m)$ be the projection
    of $k^{m+1}$ onto the first $m$ coordinates, and let $\clo{\pi(W)}$
    be the Zariski-closure
    of $\pi (W)$ in $k^m$.
    We will now examine the projection of $D$.
    Since $(f_1,\ldots, f_m)$ is algebraically independent, 
    $\pi (D)$ is Zariski-dense 
   in $k^m$, and hence $\dim (\clo{\pi (D)}) = m$.
   Since $\dim (V) \ge \dim (\clo{\pi (V)})$ 
     holds for every algebraic set $V$,
   we then obtain $\dim (\clo{D}) \ge \dim (\clo{\pi (\clo{D})}) 
                \ge \dim (\clo{\pi (D)}) = m$. 
    Seeking a contradiction, we
    suppose that  $\dim (\clo{D}) = m+1$. 

In the case  $\dim (\clo{\pi (W)}) = m$,
    we use \cite[p.\ 193, Theorem~3]{CLO:IVA}, 
     which tells $\clo{\pi (W)} = V_m ( I (W) \cap k[x_1,\ldots,x_m] )$,
    and we obtain that
    $k^m = V_m (I (W) \cap k[x_1, \ldots, x_m])$, and therefore 
    $I (W)  \cap k[x_1,\ldots,x_m] = \{0\}$. 
    Hence, \linebreak[3]
    $x_1 + I(W), \ldots, x_m + I(W)$ are algebraically independent in 
    $k[x_1,\ldots, x_{m+1}] / I(W)$.
     Since $\dim (W) \le m$, we observe that the sequence
     $(x_1 + I(W), \ldots, x_{m+1} + I(W))$ is algebraically dependent over $k$,
     and therefore, there is a polynomial $q(x_1,\ldots, x_{m+1}) \in I(W)$ with
     $\deg_{x_{m+1}} (q) > 0$. Let $r$ be the leading coefficient of 
     $q$ with respect to $x_{m+1}$, and let $(y_1, \ldots, y_m) \in k^m$
     be such that $r (y_1,\ldots, y_m) \neq 0$. Then
     there are only finitely many $z \in k$ with
     $(y_1,\ldots, y_m, z) \in W$. Since $\clo{D} = k^{m+1}$, there are then infinitely
     many $z \in k$ with  $(y_1, \ldots, y_m, z) \in D$, a contradiction to the fact that
     $g$ is $\tupBold{f}$-determined.

     In the case $\dim (\clo{\pi(W)}) \le m-1$, we
     take $(y_1,\ldots, y_{m}) \in k^{m} \setminus \pi (W)$.
     For all $z \in k$, we have $(y_1,\ldots, y_m,z) \in \clo{D}$ and
     $(y_1,\ldots, y_m, z) \not\in W$, and therefore
     all $(y_1,\ldots, y_m, z)$ are elements of $D$, a contradiction to the
     fact that $g$ is $\tupBold{f}$-determined.       
   
      Hence, we have $\dim (\clo{D}) = m$.   \qed

 \begin{thm} \label{thm:t1}
    Let $k$ be an algebraically closed field, let $\chi$ be its characteristic, 
    let $m, n  \in \N$, and
    let $(f_1,\ldots, f_m)$ be a sequence of
    polynomials in $k[\tup{t}{n}]$ that is algebraically independent over $k$.
    Then we have:
    \begin{enumerate}
        \item If $\chi = 0$, then $k \langle \tup{f}{m} \rangle  \subseteq k \FGen{f_1, \ldots, f_m} \cap k[t_1,\ldots, t_n]$.
        \item If $\chi > 0$, then $k \langle \tup{f}{m} \rangle  \subseteq \radchi (k \FGen{f_1, \ldots, f_m} \cap k[t_1,\ldots, t_n])$.
    \end{enumerate}
\end{thm}
\emph{Proof:}
    Let $g \in k \langle \tup{f}{m} \rangle$. We define
    \[ 
          D := \{ (f_1 (\vb{a}), \ldots, f_m (\vb{a}), g(\vb{a}))
                 \setsuchthat
                 \vb{a} \in k^n \},
    \] 
    we let
    $\clo{D}$ be its Zariski-closure in $k^{m+1}$, and we let 
    $W$ be an algebraic set with $\dim (W) < \dim (\clo{D})$ and 
    $\clo{D} = D \cup W$.
    By Lemma~\ref{lem:m}, we have $\dim (\clo{D}) = m$.
     Now,
     we  distinguish  cases according to the characteristic of $k$.
     Let us first suppose $\chi = 0$. 
     Let $q := \Irr (\clo{D})$ be an irreducible polynomial with
     $\clo{D} = V(q)$, and
     let $d := \deg_{x_{m+1}} (q)$.
    Since $f_1,\ldots, f_m$ are algebraically independent over $k$, 
    we have $d \ge 1$.
    We will now prove $d=1$. Suppose $d > 1$.
    We write $q = \sum_{i=0}^{d} q_i(x_1,\ldots, x_m) \, x_{m+1}^i$.       
    We recall that for a field $K$, and $f,g \in K[t]$ of positive degree, 
    the resultant 
    $\mathrm{res}_t (f, g)$ is $0$ if and only if
    $\deg (\gcd_{K[t]} (f,g)) \ge 1$ \cite[p.\ 156, Proposition~8]{CLO:IVA}. 
    Let $r := \mathrm{res}_{x_{m+1}} (q, \frac{\partial}{\partial x_{m+1}} q)$
    be the resultant of $q$ and its derivative when seen as
    elements of the ring $k(x_1,\ldots, x_{m}) [x_{m+1}]$.
    If $r = 0$, then $q$
    and $\frac{\partial}{\partial x_{m+1}} q$
    have a common divisor in
    $k (x_1,\ldots, x_m) [x_{m+1}]$ with $1 \le \deg_{x_{m+1}} (q) \le 
     d - 1$ in $k(x_1,\ldots,x_m) [x_{m+1}]$.
    Using a standard argument involving
    Gau\ss 's Lemma, we find a divisor $a$ of $q$ in
    $k[x_1,\ldots, x_{m+1}]$ such that $1 \le \deg_{x_{m+1}} (a) \le d-1$.
    This contradicts the irreducibility of $q$. 
    Hence $r \neq 0$.
 Since $\dim (\clo{\pi(W)}) \le m-1$, $r \neq 0$,
 and $q_d \neq 0$, we have
    $V(r) \cup V (q_d) \cup \pi (W) \neq k^m$.
   Thus we can choose
    $\vb{a} \in k^m$ such that $r (\vb{a}) \neq 0$, 
            $q_d (\vb{a}) \neq 0$, and
            $\vb{a} \not\in \pi (W)$.
     Let $\tilde{q} (t) := q (\vb{a}, t)$. 
    Since $\mathrm{res}_{t} (\tilde{q} (t), \tilde{q}'(t)) = r (\vb{a}) \neq 0$,
     $\tilde{q}$ has $d$ different roots in $k$, and thus
   $q (\vb{a}, x) = 0$ has $d$ distinct solutions for $x$,
    say $b_1,\ldots, b_d$.
    We will now show
     $\{ (\vb{a}, b_i) \setsuchthat i \in \{1,\ldots, d\}\} \subseteq D$.
   Let $i \in \{1,\ldots, d\}$, and suppose that $(\vb{a}, b_i) \not\in D$.
    Then 
    $(\vb{a}, b_i) \in W$, and thus $\vb{a} \in \pi (W)$,
    a contradiction. Thus all the elements $(\vb{a}, b_1), \ldots, (\vb{a}, b_d)$
    lie in $D$. 
    Since $d > 1$, this implies that $g$ is not $(f_1,\ldots,f_m)$-determined.
    Therefore we have  $d = 1$. Since $(f_1,\ldots, f_m)$ is algebraically 
    independent, the polynomial $q$ witnesses that $g$ is algebraic of degree~$1$ over
    $k \FGen{f_1,\ldots,f_m}$, and thus lies in $k \FGen{f_1,\ldots,f_m}$. 
    This concludes the case $\chi = 0$.

    Now we assume $\chi > 0$. It follows from Lemma~\ref{lem:m} 
     that for every $h \in k \langle t_1,\ldots, t_n \rangle$,
     the Zariski-closure of
     $D(h) := \{ (f_1 (\vb{a}), \ldots, f_m (\vb{a}), h(\vb{a})) \setsuchthat
                 \vb{a} \in k^n \}$ is an irreducible variety of dimension $m$
     in $k^{m+1}$. This implies that there is an irreducible polynomial
     $\Irr (\clo{D(h)}) \in k[x_1,\ldots, x_m]$ such that $\clo{D(h)} = V(\Irr(\clo{D(h)}))$.
     Furthermore, by the closure theorem \cite{CLO:IVA},  
      there is an algebraic set $W (h) \subseteq k^m$ such that
     $\dim (W (h)) \le m-1$ and $D (h) \cup W (h) = \clo{D (h)}$.  
     We will now prove the following statement by induction on 
     $\deg_{x_{m+1}} (\Irr (\clo{D(h)}))$.
     \begin{quote}
        Every $\tupBold{f}$-determined polynomial $h \in k[t_1, \ldots, t_n]$
        is an element of $\radchi (k \FGen{\tup{f}{m}} \cap k[t_1,\ldots, t_n])$.
     \end{quote}
     Let 
      \[
       d := \deg_{x_{m+1}} (\Irr (\clo{D(h)})). 
      \]
     If $d = 0$, then $f_1,\ldots, f_m$ are algebraically
     dependent, a contradiction.
     If $d = 1$, then since $f_1, \ldots, f_m$ are 
     algebraically independent, $h$ is algebraic of degree $1$ over
     $k \FGen{\tup{f}{m}}$, and thus lies in $k \FGen{\tup{f}{m}} \cap k[t_1,\ldots, t_n]$.
     Let us now consider the case $d > 1$.
     We set 
      \[ 
          e := \deg_{x_{m+1}} (\tfrac{\partial}{\partial x_{m+1}} \Irr (\clo{D(h)})).
       \]
     If $\frac{\partial}{\partial x_{m+1}} \Irr(\clo{D(h)}) = 0$, then
     there is a polynomial $p \in k[x_1,\ldots, x_{m+1}]$ such that
     $\Irr (\clo{D(h)}) = p (x_1, \ldots, x_m, x_{m+1}^{\chi})$.
         We know that $h^{\chi}$ is $\tupBold{f}$-determined, hence
     by Lemma~\ref{lem:m}, $\clo{D (h^\chi)}$ is of dimension $m$.
     Since \[p \, (f_1,\ldots, f_m, h^{\chi}) =
      \Irr (\clo{D(h)}) \, (f_1,\ldots, f_m, h) = 0, \]
     we have $p 
              \in I (D (h^\chi))$. 
     Thus $\clo{D (h^{\chi})} \subseteq V (p)$.
     Therefore, the irreducible polynomial $\Irr (\clo{D (h^{\chi})})$ divides
     $p$, and thus \[ \deg_{x_{m+1}} ( \Irr (\clo{D (h^{\chi})})) \le
     \deg_{x_{m+1}} (p) < \deg_{x_{m+1}} (\Irr (\clo{D(h)})). \] 
     By the induction hypothesis, we obtain that $h^{\chi}$ is an element
     of \linebreak[4] $\radchi (k \FGen{\tup{f}{m}} \cap k[t_1,\ldots, t_n])$.
      Therefore, $h \in \radchi (k \FGen{\tup{f}{m}} \cap k[t_1,\ldots, t_n])$.
    This concludes the case that $\frac{\partial}{\partial x_{m+1}} (\Irr (\clo{D(h)})) = 0$.
   
    If $e = 0$, we choose $\vb{a} = (a_1,\ldots, a_m) \in k^m$ such that
    \[ \frac{\partial}{\partial x_{m+1}} \Irr (\clo{D(h)}) \, \, (a_1, \ldots, a_m, 0) \neq 0,
     \] such that
    the leading coefficient of $\Irr (\clo{D (h)})$ with respect to $x_{m+1}$ does not
    vanish at $\vb{a}$, and such that
    $\vb{a} \not\in \pi (W (h))$. Then $\Irr (\clo{D(h)}) (\vb{a}, x) = 0$ has $d$ different
    solutions for $x$, say $b_1,\ldots, b_d$. Since
     $\{ (\vb{a}, b_i ) \setsuchthat i \in \{1,\ldots, d\} \} \cap W (h) = \emptyset$ because
    $\vb{a} \not\in \pi (W (h))$,
     we have $\{ (\vb{a}, b_i ) \setsuchthat i \in \{1,\ldots, d\} \} \subseteq D(h)$. Since
     $h$ is $\tupBold{f}$-determined, $d=1$, contradicting the case assumption.

     If $e > 0$, then we compute the resultant $r := \mathrm{res}^{(d,e)}_{x_{m+1}} (\Irr (\clo{D(h)}), \frac{\partial}{\partial x_{m+1}} \Irr (\clo{D(h)}))$,
     seen as polynomials of degrees $d$ and $e$ over the field $k (x_1,\ldots, x_m)$ in the variable
     $x_{m+1}$.
     As in the case $\chi = 0$, the irreducibility of $\Irr (\clo{D(h)})$ yields $r \neq 0$.
     Now we let $\vb{a} \in k^m$ be such that $r (\vb{a}) \neq 0$, the leading coefficient
     $(\Irr (\clo{D(h)}))_d$ of $\Irr (\clo{D(h)})$ with respect to $x_{m+1}$ does not vanish at $\vb{a}$, and
     $\vb{a} \not\in \pi(W(h))$. Setting $\tilde{q} (t) := \Irr(\clo{D(h)}) \, (\vb{a}, t)$, we see
     that $\mathrm{res}_t^{(d,e)} (\tilde{q} (t), \tilde{q}' (t)) \neq 0$.
     Thus $\tilde{q}$ has $d$ distinct zeroes
     $b_1,\ldots,b_d$, and then $\{ (\vb{a}, b_i) \setsuchthat i \in \{1,\ldots, d\} \} \subseteq D (h)$.
     Since $d > 1$, this contradicts the fact that $h$ is $\tupBold{f}$-determined.  \qed

\begin{thm} \label{thm:isinring0}
    Let $k$ be an algebraically closed field of characteristic~0, let $m, n \in \N$,
    and let $\tupBold{f} = (\tup{f}{m})$ be a sequence of algebraically independent
    polynomials in $k[t_1,\ldots, t_n]$.
    Then the following are equivalent:
    \begin{enumerate}
       \item \label{it:k1}
         $k \langle f_1,\ldots, f_m \rangle = k \RGen{f_1,\ldots, f_m}$.
       \item \label{it:k2}
         $\tupBold{f}$ is almost surjective.
    \end{enumerate}
\end{thm}
\emph{Proof:}
    \eqref{it:k1}$\Rightarrow$\eqref{it:k2}:
         Suppose that $\tupBold{f}$ is not almost surjective.
         Then by Lemma~\ref{lem:polies}, there are $p,q \in k[x_1,\ldots, x_m]$
         such that $p (f_1,\ldots, f_m) \mid q (f_1,\ldots, f_m)$ and
         $p \nmid q$. Let $d := \gcd (p,q)$, $p_1 := p/d$, 
         $q_1 := q/d$.
         Let $a (t_1,\ldots, t_n) \in k[t_1,\ldots, t_n]$ be such 
         that  
         \begin{equation} \label{eq:pqa} 
             p_1 (f_1,\ldots, f_m) \cdot a(t_1,\ldots, t_n) = q_1 (f_1,\ldots, f_m).
         \end{equation}
         We claim that $b(t_1,\ldots, t_n) := q_1 (f_1,\ldots, f_m) \cdot a(t_1,\ldots, t_m)$
         is $\tupBold{f}$-determined and is not an element of $k \RGen{f_1,\ldots, f_m}$.
         In order to show that $b$ is $\tupBold{f}$-determined, we let
         $\vb{c}, \vb{d} \in k^n$ such that $\tupBold{f} (\vb{c}) = \tupBold{f} (\vb{d})$.
         If $p_1 (\tupBold{f} (\vb{c})) \neq 0$, we have
         $b (\vb{c}) = q_1 (\tupBold{f} (\vb{c})) \cdot a(\vb{c}) =
                       q_1 (\tupBold{f} (\vb{c})) \cdot \frac{q_1 (\tupBold{f} (\vb{c}))}
                                                             {p_1 (\tupBold{f} (\vb{c}))} 
                  =      q_1 (\tupBold{f} (\vb{d})) \cdot \frac{q_1 (\tupBold{f} (\vb{d}))}
                                                             {p_1 (\tupBold{f} (\vb{d}))}
                   =  q_1 (\tupBold{f} (\vb{d})) \cdot a(\vb{d}) = b (\vb{d})$.
        If $p_1 (\tupBold{f} (\vb{c})) = 0$, we have
        $b (\vb{c}) = q_1 (\tupBold{f} (\vb{c})) \cdot  a(\vb{c})$.
        By~\eqref{eq:pqa}, we have $q_1 (\tupBold{f} (\vb{c})) = 0$, and thus
        $b (\vb{c}) = 0$. Similarly $b (\vb{d}) = 0$. This concludes the proof
        that $b$ is $\tupBold{f}$-determined.
        
        Let us now show that $b \not\in k \RGen {\tup{f}{m}}$.
        We have 
        \[ 
            b(t_1,\ldots, t_n)  = \frac{q_1 (f_1,\ldots, f_m)^2}{p_1 (f_1,\ldots, f_m)}.
        \]
        If $b \in k \RGen {\tup{f}{m}}$, there is $r \in k[x_1,\ldots, x_m]$ with
        $r(f_1,\ldots, f_m) = b (t_1,\ldots, t_n)$.
        Then $r(f_1,\ldots, f_m) \cdot p_1 (f_1,\ldots, f_m) = q_1 (f_1,\ldots, f_m)^2$.
        From 
        the algebraic independence
        of $(\tup{f}{m})$, we obtain 
        $r(x_1,\ldots, x_m) \cdot p_1 (x_1,\ldots, x_m) = q_1 (x_1,\ldots, x_m)^2$, hence
       $p_1 (x_1,\ldots, x_m) \mid q_1 (x_1,\ldots, x_m)^2$.
        Since $p_1, q_1$ are relatively prime, we then have $p_1 (x_1,\ldots, x_m) \mid
        q_1 (x_1,\ldots, x_m)$, contradicting the choice of $p$ and $q$.
        Hence $\tupBold{f}$ is almost surjective.

        \eqref{it:k2}$\Rightarrow$\eqref{it:k1}:
            From Theorem~\ref{thm:t1}, we obtain 
            $k \langle \tupBold{f} \rangle  \subseteq k\FGen{\tupBold{f}} \cap k[t_1,\ldots, t_n]$.
            Since $\tupBold{f}$ is almost surjective, Lemma~\ref{lem:polies} yields
            $k \FGen{\tupBold{f}} \cap k[t_1,\ldots, t_n]
                = k \RGen{\tupBold{f}}$, and thus
            $k \langle \tupBold{f} \rangle \subseteq k \RGen{\tupBold{f}}$.
            The other inclusion is obvious. \qed

\begin{thm} \label{thm:isinringchi}
    Let $k$ be an algebraically closed field of characteristic $\chi>0$, let $m, n \in \N$,
    and let $\tupBold{f} = (\tup{f}{m})$ be a sequence of algebraically independent
    polynomials in $k[t_1,\ldots, t_n]$.
    Then the following are equivalent:
    \begin{enumerate}
       \item \label{it:k1}
         $k \langle f_1,\ldots, f_m \rangle = \radchi (k \RGen{f_1,\ldots, f_m})$.
       \item \label{it:k2}
         $\tupBold{f}$ is almost surjective.
    \end{enumerate}
\end{thm}
\emph{Proof:}
     \eqref{it:k1}$\Rightarrow$\eqref{it:k2}: As in the proof of Theorem~\ref{thm:isinring0},
      we produce an $\vb{f}$-determined polynomial $b$ and relatively prime
      $p_1, q_1 \in k[x_1,\ldots, x_m]$ with $p_1 \nmid q_1$ and  
       \[
          b (t_1, \ldots, t_n) = \frac{q_1 (f_1,\ldots, f_m)^2}{p_1 (f_1, \ldots, f_m)}.
       \]
       Now suppose that there is $\nu \in \N_0$ with
       $b^{\chi^{\nu}} \in k \RGen{f_1,\ldots,f_m}$.
       Then $p_1 (f_1,\ldots, f_m)^{\chi^\nu}$  divides \linebreak[4] 
       $q_1 (f_1,\ldots, f_m)^{2 \chi^\nu}$ in
       $k \RGen{f_1,\ldots, f_m}$, and thus
       $p_1 (x_1,\ldots, x_m)$ divides $q_1 (x_1,\ldots, x_m)^{2 \chi^\nu}$ in $k[x_1,\ldots, x_m]$.
       Since $p_1$ and $q_1$ are relatively prime, we obtain
       $p_1 \mid q_1$, contradicting the choice of $p_1$ and $q_1$. 

     \eqref{it:k1}$\Rightarrow$\eqref{it:k2}:
           From Theorem~\ref{thm:t1}, we obtain 
            $k \langle \tupBold{f} \rangle  \subseteq \radchi (k\FGen{\tupBold{f}} \cap k[t_1,\ldots, t_n])$.
            Since $\tupBold{f}$ is almost surjective, Lemma~\ref{lem:polies} yields
            $k \FGen{\tupBold{f}} \cap k[t_1,\ldots, t_n]
                = k \RGen{\tupBold{f}}$, and thus
            $k \langle \tupBold{f} \rangle \subseteq \radchi (k \RGen{\tupBold{f}})$.
            The other inclusion follows from the fact that 
            the map $\varphi : k \to k$, $\varphi (y) := y^{\chi}$ is injective. \qed

\section{Function compositions that are polynomials}
 For a field $k$, let $\tupBold{f} = (\tup{f}{m}) \in (k [\tup{t}{n}])^m$, and let
 $h : k^m \to k$ be an arbitrary function.
 Then we write $h \circ \tupBold{f}$ for the function defined by
  $(h \circ \tupBold{f}) \, (\vb{a}) = 
    h(f_1 (\vb{a}), \ldots, f_m (\vb{a}))$ for all $\vb{a} \in k^n$.
 For an algebraically closed field $K$ of characteristic $\chi > 0$,
 $y \in K$ and $\nu \in \N_0$, we let $s^{(\chi^\nu)} (y)$ be the element in $K$
 with $(s^{(\chi^\nu)} (y))^{\chi^\nu} = y$; so $s^{(\chi^\nu)}$ takes
 the $\chi^{\nu}$-th root. 

\begin{thm} \label{thm:fun0}
     Let $k$ be a field, let $K$ be its algebraic
     closure, let $m, n \in \N$, 
     let $g, \tup{f}{m} \in k[t_1,\ldots, t_n]$, and let 
     $h : K^m \to K$ be an arbitrary function.
     Let $R := \tupBold{f} (K^n)$ be the range of the function from
     $K^n$ to $K^m$ that
     $\tupBold{f} = (f_1,\ldots, f_m)$ induces on $K$.
     We assume that $\dim (\clo{K^m \setminus R}) \le m-2$, and that  
      $h \circ \tupBold{f} = g$ on $K$, which means that
       \[ h (\tupBold{f} (\vb{a})) = g (\vb{a}) \text{ for all }\vb{a} \in K^n. \]
     Then we have:
     \begin{enumerate} 
        \item \label{it:c0} If $k$ is of characteristic~$0$, then
              there is $p \in k[x_1,\ldots, x_m]$ such that
              $h (\vb{b}) = p (\vb{b})$ for all $\vb{b} \in R$.
        \item \label{it:cp} If $k$ is of characteristic $\chi > 0$,
              then there is $p \in k[x_1,\ldots, x_m]$ and $\nu \in \N_0$ such that
              $h (\vb{b}) = s^{(\chi^\nu)} (p (\vb{b}))$ for all $\vb{b} \in R$.
      \end{enumerate}
\end{thm}
\emph{Proof:}
    Let us first assume that $k$ is of characteristic~$0$. 
   We observe that as a polynomial in $K[t_1,\ldots, t_n]$, 
   $g$ is $\tupBold{f}$-determined. Hence by Theorem~\ref{thm:isinring0},
   there is $q \in K[x_1,\ldots, x_m]$ such that
   $q (f_1,\ldots, f_m) = g$.
   Writing \[q = \sum_{(i_1, \ldots, i_m)\in I} \alpha_{i_1,\ldots,i_m} x_1^{i_1}\cdots x_m^{i_m}, \]
   we obtain $g = \sum_{(i_1, \ldots, i_m)\in I} \alpha_{i_1,\ldots,i_m} f_1^{i_1}\cdots f_m^{i_m}$.
   Expanding the right hand side and comparing coefficients, we see that
   $(\alpha_{i_1,\ldots,i_m})_{(\tup{i}{m}) \in I}$ is a solution of a linear
   system with coefficients in $k$. Since this system has
   a solution over $K$, it also has a solution over $k$. The solution over $k$
    provides the coefficients of a polynomial $p \in k[x_1, \ldots, x_m]$ such that
   $p (f_1,\ldots, f_m) = g$. 
   From this, we obtain that $p (f_1 (\vb{a}), \ldots, f_m (\vb{a})) = g (\vb{a})$ for all $\vb{a} \in K^n$,
   and thus $p (\vb{b}) = h (\vb{b})$ for all $\vb{b} \in R$. This completes the
   proof of item~\eqref{it:c0}.

    In the case that $k$ is of characteristic $\chi > 0$,
    Theorem~\ref{thm:isinringchi} yields a polynomial
    $q \in K[x_1,\ldots, x_m]$ and $\nu \in \N_0$ such that
   $q (f_1,\ldots, f_m) = g^{\chi^\nu}$.
  As in the previous case, we obtain
  $p \in k[x_1, \ldots, x_m]$ such that
   $p (f_1,\ldots, f_m) = g^{\chi^\nu}$. 
   Let $\vb{b} \in R$, and let $\vb{a}$ be such that
   $\tupBold{f} (\vb{a}) = \vb{b}$. Then 
    $s^{(\chi^\nu)} (p ( \vb{b} )) =
     s^{(\chi^\nu)} (p ( \tupBold{f} (\vb{a}))) 
     = g (\vb{a}) 
     = h (\tupBold{f} (\vb{a})) = h (\vb{b})$, which completes
    the proof of~\eqref{it:cp}. \qed

We will now state the special case that $k$ is algebraically
closed and $\tupBold{f}$ is surjective in the following corollary.
By a \emph{polynomial function}, we will simply mean a function induced
by a polynomial with all its coefficients in $k$.
\begin{cor} \label{cor:comp}
    Let $k$ be an algebraically closed field, let $\tupBold{f} =
    (\tup{f}{m}) \in (k[t_1,\ldots,t_n])^m$, and let
    $h : k^m \to k$ be an arbitrary function.
    We assume that $\tupBold{f}$ induces a surjective mapping from
    $k^n$ to $k^m$ and that $h \circ \tupBold{f}$ is a polynomial
    function. Then we have:
    \begin{enumerate}
      \item If $k$ is of characteristic~$0$, then $h$ is a polynomial function.
      \item If $k$ is of characteristic~$\chi > 0$, then there is $\nu \in N_0$ 
            such that $h^{\chi^\nu} : (y_1,\ldots, y_m) \mapsto
                                   h(y_1,\ldots, y_m)^{\chi^\nu}$ 
            is a polynomial function.
    \end{enumerate}
\end{cor}

\def\cprime{$'$}
\providecommand{\bysame}{\leavevmode\hbox to3em{\hrulefill}\thinspace}
\providecommand{\MR}{\relax\ifhmode\unskip\space\fi MR }
\providecommand{\MRhref}[2]{%
  \href{http://www.ams.org/mathscinet-getitem?mr=#1}{#2}
}
\providecommand{\href}[2]{#2}

\end{document}